\documentclass[24.88pt]{article}
\usepackage{latexsym,amsfonts,euscript,amssymb}

\title{Arithmetical Properties of Finite Groups\thanks{The author gratefully acknowledges the support of K.C. Wong
education foundation, Hong Kong.}}

\author{W.J. Shi\\
Department of Mathematics, Southwest China University, Chongqing 400715\\and\\
School of Mathematics, Suzhou University, Suzhou 215006, Jiangsu, P. R. China\\
E-mail: wjshi@suda.edu.cn}
\date{}

\begin{document}
\maketitle

\begin{abstract} Let $G$ be a finite group and $Ch_i(G)$ some quantitative
sets. In this paper we study the influence of $Ch_i(G)$ to the
structure of $G$. We present a survey of  author and his colleagues'
recent works.
\end{abstract}



\bigskip

Let $G$ be a finite group and $Ch(G)$ be one of the following
sets:

{\bf $Ch_1(G)=|G|$, that is, the order of $G$;

$Ch_2(G)=\pi_e(G)=\{ o(g)\,|\,\, g\in G\}$, that is, the set of
element orders of $G$;

$Ch_3(G)=cs(G)=\{ |g^G|\,\, g\in G\}$, that is, the set of
conjugacy class sizes of $G$;

$Ch_4(G)=cd(G)$, that is, the set of irreducible character degrees
of $G$.}

Our aim is to study the structure of $G$ under certain arithmetical
hypotheses of $Ch_i(G)$, $i = 1,2,3$ or $4$.

Except the above quantitative sets, we may define $Ch_5(G)$ be the
set of the maximal subgroup orders of $G$(see \cite{Li}), $Ch_6(G)$
be the set of Sylow normalizer orders of $G$(see \cite{Bi}), and the
other quantitative sets(for example, see \cite{Ab}). In this paper
we discuss the cases of  $Ch_i(G)$, $i = 1,2,3$ or $4$, especially
for the cases of $i = 1,2$.

{\bf Question  A}\,\,\em  If $Ch(G)$ is fixed, what can we say about
the structure of $G$?\em

For the set $Ch_i(G)$, $i=2,3$ or $4$, we can define a graph
$\Gamma_i(G)$ as follows: Its vertices are the primes dividing the
numbers in $Ch_i(G)$; and two distinct vertices $p,q$ are connected
if $pq\,|\,m$ holds for $m \in Ch_i(G)$.

{\bf Question B}\,\,\em If we know the information of graph
$\Gamma_i(G)$, $i=2,3$ or $4$, what can we say about the structure
of $G$?\em

In our characterization using the element orders, the graph
$\Gamma_2(G)$(prime graph) and Gruenberg-Kegel theorem on groups
with disconnected prime graphs(see \cite{Will}) play an important
role.

For any $a \in Ch_i(G)$, $i=2,3$ or $4$, for example, $a \in
Ch_2(G)$, we define $M_2(a)$, the multiplicity of $a$ in $G$ as the
number of elements of order $a$. Also, $DC_{3,2}(G)=
|Ch_3'(G)|-|Ch_2(G)|$, the difference of conjugacy classes number
$|Ch_3'(G)|$ and same order classes number in $G$. Furthermore,
$QC_{1,3}(G)= |G|/|Ch_3'(G)|$, the quotient of $|G|$ and the number
of conjugacy classes.

{\bf Question C}\,\,\em  If $Ch(G)$ and $M(a)$, for all $a \in G$,
are known what can we say about the structure of $G$? If $DC(G)$ is
"small", what can we say about the structure of $G$? If $QC(G)$ is
known what can we say about the structure of $G$?\em

For example, $DC_{3,2}(G)=0$, that is, any elements of a finite
group $G$ with same order are conjugate. This is Syskin's problem,
and some group theorists have proved that $G\cong 1, Z_2$, or $S_3$
using the classification theorem of finite simple
groups(see\cite{FS}, \cite{Zh}). In \cite{DS2} we classified all
finite groups of $DC_{3,2}(G)=1$, and studied such finite groups in
which elements of the same order outside the center are
conjugate(see \cite{QSY}).

{\bf Question D}\,\,\em Let $b_i(G)=max \{Ch_i(G)\}$, $i=2,3$ or
$4$. If some information of $b_i(G)$ is known, what can we say about
the structure of $G$?\em

\bigskip
For Problem A we know there exists some very famous results for
$Ch_1(G)=|G|$, for example, {\bf Sylow's theorem, the odd order
theorem}, and {\bf Burnside's $p^aq^b$ theorem} etc.(see\cite{B},
\cite{FT} and some important bibliography listed in \cite{G}).

{\bf Burnside's $p^aq^b$ theorem} implies that if $G$ is a
non-Abelian simple group, then $|\pi(G)|\geq 3$.

A finite group $G$ is called a $K_n$-group if $|\pi(G)|= n$. There
are only eight simple $K_3$-groups(see\cite{He}). In \cite{Sh1} we
classified all simple $K_4$-groups, but we do not know whether the
number of $K_4$-groups is finite or not. This problem depends on the
solutions of some special Diophantine equations(also see\cite{BCM}).

Twenty years ago, we studied such finite groups in which every
non-identity element has prime order, i.e. finite EPO groups and got
an interesting result: $G\cong A_5$ if and only if
$\pi_{e}(G)=\{1,2,3,5\}$(That is, we may characterize the integral
property only using the local property, see \cite{SY}).

After this, we developed such characterizations for all finite
simple groups using the "two orders", and for the finite nonsolvable
groups only using the "set of elements orders"(or "spectrum"). That
is, we researched the following characterizations of two kinds.

(1) Characterizing all finite simple groups unitization using only
the two sets $Ch_1(G)$ and $Ch_2(G)$.

We have finished the above works except $B_n$, $C_n$ and $D_n$ ($n$
even)(see\cite{Sh2}, \cite{Sh3}, \cite{SB1}, \cite{SB2}, \cite{SB3},
\cite{CS} and \cite{XS}).

(2) Many finite simple groups are characterizable using only the set
$Ch_2(G)$.

The most recent version of the latter characterization is presented
in Table 1 of \cite{M3}.

Let $|G|=n$ and $f(n)$ denote the numbers of $G$, pairwise
non-isomorphic, such that $|G|=n$. It is easy to see that $f(n)=1$
if and only if $(n,\varphi(n))=1$, where $\varphi(n)$ is a Euler
function of $n$. Also the solutions for $f(n)=2,3,$ and $4$ are
found(see\cite{W}).

The following question is posed: For any integer $k$, is there a
solution for $f(n)= k$?

Now we consider substituting $|G|$ by the set $\pi_{e}(G)$.
$\pi_{e}(G)$ is a set of some positive integers. Similarly, for a
set $\Gamma$ of positive integers, let $h(\Gamma)$ be the number of
isomorphic classes of finite group $G$ such that
$\pi_{e}(G)=\Gamma$.

If $\Gamma=\pi_{e}(G)$, then we have $h(\pi_{e}(G))\geq 1$. If
$\Gamma=\{1,2,3,5\}$, then $h(\Gamma)=1$. Conversely, which groups
$G$ satisfy $h(\pi_{e}(G))=1$? Such groups are called {\bf
characterizable groups} or {\bf recognizable groups}.

With respect to characterizable groups, summarizing many scholars'
works, we have the following results:

{\bf Theorem 1}\,\,\em  The following groups are characterizable
groups:

(a) Alternating groups $A_n$, where $n=5,16,p,p+1,p+2,$ and $p\geq
7$ is a prime; Symmetric group $S_n$, where
$n=7,9,11,12,13,14,19,20,23,24$.

(b) Simple groups of Lie type $L_2(q), q\neq 9$, $L_3(2^m)$,
$L_n(2^m)$, $n=2^k>8$, $F_4(2^m)$, $U_3(2^m)$,$m\geq 2$, series of
simple groups of Suzuki-Ree type $Sz(2^{2m+1})$,
${}^2G_2(3^{2m+1})$, ${}^2F_4(2^{2m+1})$; $S_4(3^{2m+1}) (m>0)$,
$G_2(3^m)$; $L_3(7)$, $L_4(3)$, $L_5(3)$, $L_5(2)$, $L_6(2)$,
$L_7(2)$, $L_8(2)$, $U_3(9)$, $U_3(11)$, $U_4(3)$, $U_6(2)$,
$S_6(3)$, $O^{-}_8(2)$, $O^{-}_{10}(2)$, ${}^3D_4(2)$, $G_2(4)$,
$G_2(5)$, $F_4(2)$, ${}^2F_4(2)'$, ${}^2E_6(2)$; and nonsolvable
groups $PGL_2(p^m)$, $m>1$, $p^m\neq 9$, $L_2(9).2_3$ $(\cong
M_{10})$, $L_3(4).2_1$.

(c) All sporadic simple groups except $J_2$.\em

\cite{KM} extended the characterization of the above some groups
from finite groups to periodic groups.

For the case of $h(\pi_{e}(G))= \infty$ we have

{\bf Theorem 2}\,\,\em  If all the minimal normal subgroups of $G$
are elementary(especially $G$ is solvable), or $G$ is one of the
following: $A_6$, $A_{10}$, $L_3(3)$, $U_3(3)$, $U_3(5)$, $U_3(7)$,
$U_4(2)$, $U_5(2)$; $J_2$; $S_4(q) (q\neq 3^{2m+1}$ and $m>0)$, then
$h(\pi_e(G))=\infty.$\em

In the case of $k$-recognized groups, we have found infinite pairs
of 2-recognizable groups as follows:

(a) $L_3(q)$, $L_3(q)\langle\theta\rangle$, where $\theta$ is a
graph automorphism of $L_3(q)$ of order 2, $q=5,29,41$, or
$q\equiv\pm2(mod\ 5)$ and $(6, (q-1)/2)=2$(see \cite{M1}, \cite{MD}
and some references listed in \cite{MD});

(b) $L_3(9)$, $L_3(9).2_1$(see \cite{CS1});

(c) $S_6(2)$, $O_8^+(2)$(see \cite{M2} and \cite{ST});

(d) $O_7(3)$, $O_8^+(3)$(see \cite{ST});

(e) $L_6(3)$, $L_6(3)\langle\theta\rangle$,where $\theta$ is a graph
automorphism of $L_6(3)$ of order 2(see \cite {Va});

(f) $U_4(5)$, $U_4(5)\langle\gamma\rangle$,where $\gamma$ is a graph
automorphism of $U_4(5)$ of order 2(see \cite {Va}).

For any $r>0$, $h(\pi_e(L_3(7^{3^r})))=r+1$, and these $r+1$ groups
are $L_3(7^{3^r})\langle\rho\rangle$, where $\rho$ is a field
automorphism of $L_3(7^{3^r})$, $k=0,1,2,\ldots,r$(see \cite{Za}).

{\bf Problem 1}\,\,\em For the cases of $B_n(q)$ and $C_n(q)$, $q$
odd, we have $|B_n(q)|$ = $|C_n(q)|$, and $B_n(q)$ is not isomorphic
to $C_n(q)$. How to distinguish them using the set $\pi_e(G)$? \em

{\bf Problem 2}\,\,\em Find new characterizable simple groups within
a particular class of finite simple groups.\em

Considering the independent number of the prime graph A.V. Vasil'ev
and E.P. Vdovin find recently a new approach which makes possible to
study the case of a finite simple group with the connected prime
graph(see \cite{VV}, also see \cite{VG} and \cite{GSV} for its
application).

{\bf Problem 3}\,\,\em Whether or not there exist two section-free
finite groups $G_1$ and $G_2$ such that $\pi_e(G_1)=\pi_e(G_2)$ and
$h(\pi_e(G_1))$ is finite?\em

{\bf Problem 4}\,\,\em For a enough large positive integer $n$, is
alternating group $A_n$ characterizable?\em

{\bf Problem 5}\,\,\em For any $n\geq 3$, does $h(L_n(2))=1$(see
\cite{GLMMV})?\em

In 1987, when the author communicated the first
characterization(that is, Characterizing all finite simple groups
unitization using only the two sets $|G|$ and $\pi_e(G)$) with Prof.
J.G. Thompson, he put forward the following questions and
conjectures in his letters:

For any finite group $G$ and any integer $d>0$, let $G(d) = \{x \in
G; x^d = 1\}$. Two finite groups $G_1$ and $G_2$ are of {\bf the
same order type} if and only if $|G_1(d)| = |G_2(d)|$, $d = 1, 2,
\cdots$(That is, $Ch_2(G_1)= Ch_2(G_2)$ and $M_2(G_1)= M_2(G_2)$).

{\bf Problem 6 (J.G. Thompson)}\,\,\em Suppose $G_1$ and $G_2$ are
groups of the same order type. If  $G_1$ is solvable, is  $G_2$
necessarily solvable?\em

The problem that the solvability of groups in which the number of
elements whose orders are largest are given, induced by Thompson's
problem, interested many Chinese group-theory specialists. They
proved the following results:

{\bf Theorem 4}\,\,\em Let $G$ be a finite group and $b'_2(G)$ be
the number of elements of maximal order in $G$. If $b'_2(G)$ = odd,
$32, 2p, 4p, 6p, 8p, 10p, 2p^2, 2p^3, 2pq$ ($p, q$ are primes), or
$b'_2(G) = \varphi(k)$, where $\varphi(k)$ is the Euler function of
maximal order $k$, then $G$ is solvable except $G\cong S_5$.\em

{\bf Corollary}\,\,\em  For the above cases, Thompson's problem is
affirmative.\em

{\bf Problem 7 (J.G. Thompson's conjecture)}\,\,\em  Let $G$, $H$ be
two finite groups with $Ch_3(G)= Ch_3(H)$ ($cs(G)=cs(H)$). If $H$ is
a non-Abelian simple group and $Z(G)$ = 1, then $G \cong H$.\em

G.Y. Chen has introduced the concept of {\bf order components} and
proved that Thompson's conjecture holds if $H$ is a finite simple
group at least three prime graph components(see \cite{C}).

Now we consider the case of "small" difference number $DC(G)$.

{\bf Theorem 5}\,\,\em Let $G$ be a finite group. Then
$DC_{3,2}(G)=1$(i.e., $G$ have one and only one same order class
containing two conjugacy classes of $G$) if and only if $G\cong A_5,
L_2(7), S_5$, $\,S_4,\,A_4,\,D_{10},\,Z_3,\,Z_4,\,Hol(Z_5)$ or
$[Z_3]Z_4 $(see \cite{DS2}).\em

{\bf Theorem 6}\,\,\em Let $G$ be a finite group. \\
(1) If $G$ is non-Abelian, then $QC_{1,3}(G) \geq 8/5$, and
$QC_{1,3}(G)= 8/5$ if and only if $G = P \times A$, where $A$ is
abelian with odd order and $P$ is a specific non-Abelian 2-group.\\
(2) If $Z(G) = 1$, then $QC_{1,3}(G) \geq 2$, and $QC_{1,3}(G) = 2$
if and only if $G \cong S_3$.\\
(3) If $G$ is non-Abelian simple, then $QC_{1,3}(G) \geq 12$, and
$QC_{1,3}(G) = 12$ if and only if $G \cong A_5$(see \cite{SX}).\em

Some papers improved the above result(see \cite{Zho} and \cite{Du}).

\bigskip
For the case of $Ch_3(G)$ and $Ch_4(G)$, we may pose the similar
problems. For example, {\bf which positive integers set can become
$Ch_3(G)$ or $Ch_4(G)$ for some groups $G$?} Also, we may define the
corresponding graphs. In \cite{Lew}, M.L. Lewis presented an
integral overview. The following results are just related with the
author's joint works.

{\bf Theorem 7} \,\, (a) (\cite{Ito1})\em If $G$ is a finite group
with $|cs(G)|= 2$, then (1) $G=P\times A$ with $P$ a $p$-group and
$N$ Abelian. (2) $P/Z(P)$ has exponent $p$. Also, if $dl(P)\leq 2$,
then $c(P)\leq 3$.\em

(b) (\cite{BMHQS}) \em Let $G$ be a finite group. Then $|cd(G)|= 2$
if and only if one of the following is true: (1) $G$ possesses a
normal and Abelian subgroup $N$ with $|G:N|=q$, $q$ is a prime. (2)
$G/Z(G)$ is a Frobnius group with kernel $(G'\times Z(G))/Z(G)$ and
a cyclic complement. (3) $G = P\times A$, where  $A$ is an Abelian
and $P$ is a $p$-group with $|cd(P)|=2$. \em

For the case (a) it is proved that the nilpotent class of $G$ is at
most 3(see \cite{Ishi}).

{\bf Theorem 8}\,\,  (a) (\cite{Ito2, Ito3, Ito4}) \em If
$|cs(G)|\leq 3$, then $G$ is solvable; If $G$ is simple with
$|cs(G)|=4$, then $G\cong PSL(2, 2^m)$ (and conversely); If $G$ is
simple with $|cs(G)|=5$, then $G\cong PSL(2,q)$, where $q$ is an odd
prime power greater than 5(and conversely).\em

(b) (\cite{Isa})\em If $|cd(G)|=k\leq 3$, then $G$ is solvable. If
$G$ is solvable and $|cd(G)|=k\leq 4$, then $dl(G)\leq k$.\em

(c) (\cite {QS})\em If $cd(G)=\{1,m,n, mn\}$, then $G$ is solvable
and one of the following is true: (1) $dl(G)\leq 3$. (2)
$cd(G)=\{1,3,13,39\}$. (3) $cd(G)=\{1, p^a, p^b, p^{a+b}\}$.\em

In \cite {MM}, the authors classified nonsolvable groups with four
irreducible character degrees. Furthermore, it is proved that if
$cs(G)= \{1, m, n, mn\}$ and $(m, n)=1$ then $G$ is solvable(see
\cite {BF}).\em

{\bf Acknowledegement.} The author would like to thank Dr. G.H.
Qian, Dr. A. Moret\'o and the referee for their many help.

\end{document}